\documentclass[a4paper]{amsart}

\usepackage{amssymb}
\usepackage{amsmath,amsthm}
\usepackage{cite}

\usepackage{enumerate,paralist}
\usepackage{amstext}
\usepackage{psfrag}
\usepackage{dsfont}
\usepackage[dvips]{epsfig}
\usepackage[english]{babel}
\usepackage{url}
\usepackage{graphicx}
\usepackage{dcolumn}
\usepackage{bm}
\usepackage{float}
\usepackage{hyperref}
\usepackage{color}
\usepackage{epstopdf}
\usepackage{cleveref}
\usepackage[svgnames]{xcolor}
\hypersetup{hidelinks,colorlinks=true,allcolors=DarkBlue}

\theoremstyle{plain}
\newtheorem{theorem}{Theorem}
\newtheorem{corollary}{Corollary}

\theoremstyle{definition}
\newtheorem{definition}{Definition}

\DeclareMathOperator{\sign}{sign}

\def\text#1{\mbox{\rm #1\,}}

\begin{document}
\frenchspacing

\title[Motion of a system of oscillators]
{Motion of a system of oscillators under the generalized dry friction control}

\author{Alexander Ovseevich}
\address
{
Institute for Problems in Mechanics, Russian Academy of Sciences \\
119526, Vernadsky av., 101/1, Moscow, Russia
}
\email{ovseev@ipmnet.ru}

\author{Aleksey Fedorov}
\address
{
Institute for Problems in Mechanics, Russian Academy of Sciences \\
119526, Vernadsky av., 101/1, Moscow, Russia
\&
Bauman Moscow State Technical University,
105005, 2nd Baumanskaya str., 5, Moscow, Russia
}
 \email{akf@rqc.ru}

\maketitle

\begin{flushright}
    {\it To Boris T. Polyak \\ on occasion of his jubilee}
\end{flushright}

\begin{abstract}
The problem of the existence and uniqueness of the motion of the system of an arbitrary number linear oscillators under a generalized dry-friction type control is studied. 
This type of control arises in the problem of steering the system to equilibrium. 
The problem of existence and uniqueness of motion under the suggested control is resolved within the framework of the DiPerna–Lions theory of singular ordinary differential equations.


\end{abstract}

\section{Introduction}

The content of the paper is based on the author's talk at XII All-Russian conference on control problems \cite{doklad}.
B.T.~Polyak have expressed a great interest to the author's presentation.
This is why we are especially pleased to publish this paper in the jubilee issue of the ``Automation and Remote Control''.

The control of oscillatory systems is one of the most important issues of the optimal control theory \cite{akulenko}.
A classical achievement in this area is the explicit construction of the feedback control for the minimum-time damping of a single linear oscillator \cite{pont}.
From the point of view of the canonical system of the maximum principle, the problem is completely integrable.

The subject of the present paper is a more complicated problem of control for a system of an arbitrary number $N$ of linear oscillators with different eigenfrequencies $\omega_i$,
\begin{equation}\label{syst0}
\begin{array}{l}
    \dot{x_i}=y_i \\
    \dot{y_i}=-\omega_i^{2}{x_i}+{u}, \quad |u|\leq1, \quad i=1,\dots,N,
\end{array}
\end{equation}
which is probably  not completely integrable.

Authors do not expect to design an optimal control analytically, however the problem of damping of oscillations of system (\ref{syst0}) is to be solved.
A standard way to damp oscillation goes via friction.
The straightforward application of the Coulomb law leads to the vector control $u_i=-\sign{y_i}$.
The {\it scalar} control we are using is a generalization of the dry friction in the sense that it has the form
\begin{equation}
    u=-\sign{\sum_{i=1}^{N}\lambda_iy_i},
\end{equation}
where $\lambda_i$ are some coefficients.
This type of control arises naturally in author's papers \cite{ovseev0,ovseev2,ovseev3} based on the study of asymptotic behavior of reachable sets of the system (\ref{syst0}).

It is well known that the use of dry friction, although it helps damping oscillations, does not necessarily leads to a complete stop of the system:
there arise standstill zones where the system is not moving at all, although the terminal set is not yet reached.
In \cite{ovseev0,ovseev2,ovseev3} an approach to design of quasi-optimal control based on a combination of three different strategies has been suggested.
At low energies, 
the problem is resolved by using the common Lyapunov functions \cite{ovseev3}, whereas at high and intermediate energies the  control used takes the form of the generalized dry friction.

A characteristic feature of control problems is the discontinuity of the right-hand sides of the differential equations of motion.
One encounters the same phenomenon in the study of the motion under the dry friction.
The question of existence of the motion is traditionally resolved by using the Filippov theory of existence of  solutions of differential inclusions \cite{filipp}.
However, the problem of uniqueness of solutions goes beyond the Filippov theory, although the intuitive concept of motion stipulates an unambiguous determination of the trajectories of the system by the control law.

In this paper, the problem of existence and uniqueness of the motion of the system of oscillators under the action of a control in the dry friction from is investigated.
We explain how the control in the dry friction form arises from the structure of asymptotic behavior of reachable sets.
In the present work, the central issue is the existence of uniquely defined motion under the control.
This problem is resolved in the framework of the DiPerna--Lions theory of singular ordinary differential equations (ODE) \cite{diperna, ovseev_diperna}.
In this paper, we confine ourselves to the statement of the existence and uniqueness theorems.
Proofs can be found in \cite{ovseev0}.

Consider a system of an arbitrary number of linear oscillators under a common bounded scalar control.
The systems has the from:
\begin{equation}\label{syst1}
    \dot{x}={A}x+{B}u,\qquad    x\in{\mathbb V}={\mathbb R}^{2N},\qquad u\in{\mathbb U}={\mathbb R},\qquad |u|\leq1,
\end{equation}
where the matrix ${A}$ and the vector ${B}$ are as follows:
\begin{equation}\label{syst2}
    {A_i}=
    \left({\begin{array}{*{20}c}
    0   &   1  \\
    { - \omega_i^{2}} & 0 \\
    \end{array}}\right), \qquad
    {A}={\rm diag}(A_i), \qquad
    {B_i} = \left( \begin{gathered}
    0 \hfill \\
    1 \hfill \\
    \end{gathered}\right), \qquad
    {B}=\oplus B_i.
\end{equation}

It follows from the Bellman dynamic programming principle  that the minimum-time control has the from \cite{bellman}:
\begin{equation}\label{control_o}
    u(x)=-\sign\langle B,p(x)\rangle, \qquad p=p(x)=\frac{\partial T}{\partial x}(x).
\end{equation}
Here $T(x)$ the optimal damping time, $p$ is the outer normal to the reachable set from zero $\mathcal{D}(T(x))$,
whose boundary passes through $x$, and the angle brackets stand for standard scalar product in $\mathbb{R}^{2N}$.
Recall that the reachable set $\mathcal{D}(T)$ is defined as the set of point reachable from zero at the time $T$ along the trajectories of system (\ref{syst1}).

In the real life, the reachable sets are quite complicated.
This is why we will design a control by using Eq. (\ref{control_o}), where $p$ is the normal to an {\it approximate} reachable set such that its boundary passes through $x$.
According to the asymptotic theory of reachable sets, for system (\ref{syst1}) a good approximation to $\mathcal{D}(T)$ as $T\to\infty$ is provided by a set of the form $T\Omega$,
where $\Omega$ is a time-independent convex body \cite{ovseev4,ovseev5,ovseev6}.
As it is well known, a closed convex body $M$ is uniquely determined by its support function \cite{schneider},
\begin{equation}
    {H}_{M}(\xi)=\sup_{x\in M}\langle{\xi, x\rangle}.
\end{equation}

The language of support functions \cite{ovseev4} is convenient for a precise statement on the asymptotic behavior of reachable sets:
\begin{theorem}\label{support0}
    Suppose that the momentum $p$ is written in the form $p=(p_i)$,
    where $p_i=(\xi_i,\eta_i)$, ${i=1,\dots,N}$,
    $\xi_i$ is the dual variable for $x_i$,
    $\eta_i$ is the dual variable for $y_i$, and
    $z_i=(\eta_i^2+{\omega_i^{-2}}{\xi_i^{2}})^{1/2}$.
    In the absence of resonances, {\it i.e.}, non-trivial relations between eigenfrequencies of from
    \begin{equation}\label{reson}
        \sum_{i=1}^Nm_i\omega_i=0, \mbox{ ãäå } 0\neq m=(m_1,\dots,m_N)\in {\mathbb{Z}}^N,
    \end{equation}
    the support function $H_T$ of the reachable set $\mathcal{D}(T)$ has as $T\to\infty$ the asymptotic form
    \begin{equation}\label{approxN0}
        {H}_T(p)=\frac{T}{(2\pi)^N}\int_{0}^{2\pi}\dots\int_{0}^{2\pi}\left|\sum_{i=1}^N z_{i}\cos\varphi_{i}\right|d\varphi_{1}\dots d\varphi_{N}+o(T)=T\mathfrak{H}(z)+o(T),
    \end{equation}
    and the support function of the closed compact $\Omega$ is given by the main term $\mathfrak{H}(z)$.
\end{theorem}

\subsection{Formula for the support function}

The support function ${H}_\Omega(p)$ of the convex body $\Omega$ is given by the main term of asymptotic equality (\ref{approxN0}):
\begin{equation}\label{approx2N}
    {H}_\Omega(p)=\mathfrak{H}(z)=\int\left|\sum_{i=1}^N z_{i}\cos\varphi_{i}\right|d\varphi, \mbox{ ãäå } z=(z_1,\dots,z_N)\in{\mathbb R}^N.
\end{equation}

If $N=1$, we obtain $\mathfrak{H}(z)=\frac2\pi|z|$.
In the case of two oscillators $N=2$, the function
\begin{equation}
    \mathfrak{H}(z)=\int\left|z_{1}\cos\varphi_{1}+z_{2}\cos\varphi_{2}\right|d\varphi
\end{equation}
can be expressed via elliptic integrals  \cite{ovseev0}.
In general, function
\begin{equation}
    \mathfrak{H}(z)=\frac{1}{(2\pi)^N}\int\limits_{\{|t_i|\leq1\}}{\left|\sum_{i=1}^N z_{i}t_{i}\right|}{\prod_{i=1}^N(1-t_i^2)^{-1/2}}dt_1\dots dt_N
\end{equation}
is an Euler-type integral that defines a generalized hypergeometric function in the sense of I.M.~Gelfand.

\subsection{Equation defining the control}
The vector $p$ is normal to the boundary  $\partial\Omega$ at the point $\frac{\partial H_\Omega}{\partial p}(p)$.
Therefore, the normal vector to the approximate reachable set $\rho\Omega$ such that its boundary passes through $x$ can be defined via the equation
\begin{equation}\label{approx30}
    \rho^{-1}x=\frac{\partial {H}_\Omega(p)}{\partial p}=\frac{\partial \mathfrak{H}(z)}{\partial z}\frac{\partial z}{\partial p},
\end{equation}
where $p\in\mathbb{R}^{2N}$ and $\rho>0$ are unknown.
We note that the support function ${H}_{\Omega}$ is differentiable, and equation (\ref{approx30}) has only one solution due to the smoothness of the boundary of the set $\Omega$ established in \cite{ovseev5}.

\subsubsection{Duality}

More generally, the relation between the support function $H=H_\Omega(p)$ and the function $\rho=\rho(x)$ in equation (\ref{approx30}) is similar to the Legendre transform:
\begin{equation}\label{Legendre}
    \langle x,p\rangle=\rho(x) H(p),\quad \rho(x)=\max_{H(p)\leq1}\langle x,p\rangle,\quad H(p)=\max_{\rho(x)\leq1}\langle{x,p}\rangle,
\end{equation}
and the corresponding point transformation $x\rightleftarrows p$ has the form
\begin{equation}\label{Legendre2}
    x=\rho(x)\frac{\partial H}{\partial p}(p),\qquad p=H(p)\frac{\partial \rho}{\partial x}(x).
\end{equation}
Here $p$ and $x$ are maximum points in equations (\ref{Legendre}).
These relations make sense if the functions $H$ and $\rho$ are norms of their arguments.
The differentiation in (\ref{Legendre2}) refers to taking subgradients.
If one of the function $H$ or $\rho$ is strictly convex and differentiable outside zero, then the other is so.
Then relations (\ref{Legendre2}) hold in the classical sense.

By using the function $\rho=\rho(x)$, one can rewrite control (\ref{control_o}) in the form
\begin{equation}\label{approx4}
    u(x)=-\sign\left\langle B,\frac{\partial\rho}{\partial x}\right\rangle.
\end{equation}
Equation (\ref{approx4}) is used even in the resonant case, when the asymptotic equality (\ref{approxN0}) does not work.
Note that $u(x)$ is a multivalued function, because $\sign(0)$  can take any value in the interval $[-1,1]$.

\subsubsection{Hamiltonian structure}

For the considered control, the function $\rho=\rho(x)$ from equation (\ref{approx30}) plays the same role as the optimal damping time $T(x)$ for the optimal control.
The ``canonical'' momentum in (\ref{approx30}) is $p=-{\partial\rho}/{\partial x}$.
The function $\rho(x)$ is the norm of the vector $x$ in the metric such that the body $\Omega$ is the unit ball wrt to it.
This is a smooth function outside zero.

Furthermore, a ``maximum principle'' holds: the compound vector $(x,p)$, where $p={\partial\rho}/{\partial x}$, satisfies the {\it Hamiltonian} system
\begin{equation}\label{attractor_syst22p110}
    \dot x=Ax+ Bu,\quad u=\sign\langle{B,p}\rangle,
    \quad \dot p=-A^*p+\frac{\partial^2\rho}{\partial x^2}B\,\sign\left\langle{B,\frac{\partial\rho}{\partial x}}\right\rangle,
\end{equation}
where the Hamiltonian has the form
\begin{equation}\label{hamilton2}
    \mathcal{H}=\langle{Ax,p}\rangle+|\langle{B,p}\rangle|-\left|\left\langle{B,\frac{\partial\rho}{\partial x}}\right\rangle\right|.
\end{equation}
The Hamiltonian vanishes along  trajectories of the controlled motion in  complete similarity with the case of the optimal control.

\section{Motion under the control}

Formally speaking, the motion of the oscillator system under control (\ref{approx4}) is described by the differential equation
\begin{equation}\label{sing_sys}
    \dot x=Ax- B\,\sign\left\langle{B,\frac{\partial\rho}{\partial x}}\right\rangle.
\end{equation}
In fact, as it is already mentioned, the right-hand side of this equation is not defined uniquely ($\sign(0)=[-1,1]$).
Therefore, we are dealing with a differential inclusion.

Fillippov theory \cite{filipp} says that the Cauchy problem for differential inclusion (\ref{sing_sys}) is solvable for any initial condition $x(0)$,
{\it i.e.} there is a function of time $x(t)$, which is absolutely continuous and satisfies (\ref{sing_sys}) at points of differentiability.
This follows from the fact that the right-hand side $f(x)$ is, first, linearly bounded $|f(x)|\leq C|x|$.
Second, it takes convex and compact values.
Finally, the function $f(x)$ is continuous as a multivalued function: if $y_n\in f(x_n)$ and $x_n\to x$, the $y\in f(x)$, where $y$ is any limit point of the sequence $y_n$.
Generally speaking, the Filippov ``motion'' is not defined uniquely: several trajectories from a single starting point are possible.

The central result of the paper is that the motion under control (\ref{approx4}) can be defined uniquely.
Toward this end, the DiPerna--Lions theory of singular ODE \cite{diperna} is used.
The general idea is to build a global phase flow, perhaps, not everywhere defined, instead of finding an individual solution of the Cauchy problem for each initial condition.

\subsection{DiPerna--Lions theory}

If $b(x)$ is a Lipschitz function, then the Cauchy problem for ODE
\begin{equation}\label{ODE}
    \dot x=b(x),\quad x(0)=x_0
\end{equation}
and for the partial differential equation
\begin{equation}\label{PDE}
    \frac{\partial v}{\partial t}=\sum b_i(x)\frac{\partial v}{\partial x_i},
    \quad
    v(x,0)=v_0(x)
\end{equation}
are equivalent.
Using the method of characteristics, one can show that the solution $v$ of (\ref{PDE}) is given by the formula
\begin{equation}\label{flow}
    v(x,t)=v_0(\phi_t(x)),
\end{equation}
where $\phi_t$ is the phase flow for (\ref{ODE}).

In the classical work of DiPerna and Lions \cite{diperna}, the Lipschitz condition on $b$ is substantially relaxed, and it is shown that the solution of problem (\ref{PDE}) still exists and is unique and is given by formula (\ref{flow}),
where $\phi_t$ is the measurable flow.
Roughly speaking, it is shown that it suffices that the Lipschitz condition holds in the integral sense, not pointwise.

\begin{definition}
    A weak bounded solution $v$ of the Cauchy problem is called renormalized solution,
    if for any smooth function $\beta:\mathbb{R}\to\mathbb{R}$ function $\beta(v)$ is again a weak solution.
\end{definition}

Using this definition, we formulate the following theorem:
\begin{theorem}\label{hamiltonian_motion}
    Suppose that  (extended) DiPerna--Lions conditions hold in the Cauchy problem for the transport equation (\ref{PDE}):
    \begin{equation}\label{dip_le2a2}
        {\rm div} b\in L^\infty, \quad  b\in {W_{ *\,{\rm  loc}}^{1,1}}=BV_{\rm  loc}, \quad \frac{b(x)}{1+|x|}\in L^\infty+L^1,
    \end{equation}
    where $BV_{\rm  loc}={W_{*\,{\rm loc}}^{1,1}}$ is the space of functions such that their derivatives are locally finite
    measures.
    Then there exists a unique measurable flow $\phi_t(x)$ such that if $v_0(x)$ is a bounded measurable function, then the function
        \begin{equation}
        v(x,t)=v(\phi_t(x))
    \end{equation}
    is the unique renormalized solution of the Cauchy problem (\ref{PDE}).
\end{theorem}
We note that the Lipschitz condition can be restated as ${\partial b}/{\partial x}\in L^\infty$.
This theorem can be easily proven by methods of \cite{ovseev_diperna}, where a strengthening of the results of \cite{diperna} is suggested.

\begin{corollary}\label{ham_setup}
    The Cauchy problem for the transport equation that corresponds to Hamiltonian system (\ref{attractor_syst22p110}) and a bounded initial condition $v_0(x,p)$ has a unique renormalized solution $v$.
    The solution has the form
    \begin{equation}
        v(x,p,t)=v_0(\phi_t(x,p)),
    \end{equation}
    where
    \begin{equation}
        \phi_t:{\mathbb{R}}^{4N}\to{\mathbb{R}}^{4N},
        \quad t\in{\mathbb{R}}
    \end{equation}
    is a uniquely defined measurable flow.
    Each curve
        \begin{equation}
        t\mapsto (x(t),p(t))=\phi_t(x,p)
        \end{equation}
    is absolutely continuous and satisfies (\ref{attractor_syst22p110}).
\end{corollary}
Theorem \ref{hamiltonian_motion} and its corollary are useful, but they do not define any flow in the phase space
$\mathbb{R}^{2N}$ of system (\ref{sing_sys}), because in the extended symplectic space $\mathbb{R}^{4N}$ the phase space
has measure zero.

\subsection{Existence and uniqueness of the motion}
The Cauchy problem for the transport equation related to ODE (\ref{sing_sys}) is as follows:
\begin{equation}\label{transport}
\frac{\partial v}{\partial t}=\left\langle Ax-B\sign\left\langle{B,\frac{\partial\rho}{\partial
x}(x)}\right\rangle,\frac{\partial v}{\partial x}\right\rangle,\quad v(x,0)=v(x).
\end{equation}
Our main result claims that in the phase space ${\mathbb{R}}^{2N}$ of the system (\ref{sing_sys}) we can define a
\textit{semiflow} which is continuous, uniquely defined everywhere, and it is related to the transport equation
(\ref{transport}) in a way the flow from Theorem \ref{hamiltonian_motion} is related to (\ref{PDE}):
\begin{theorem}\label{motion2}
    There exists a unique continuous semiflow $\phi_t(x),\,t\geq0$ such that $v(x,t)=v(\phi_t(x))$ is the unique renormalized
    solution of the Cauchy problem for the transport equation (\ref{transport}).
    Every curve $x(t)=\phi_t(x)$ is absolutely continuous and differential inclusion (\ref{sing_sys}) holds.
\end{theorem}
The proof is rather lengthy and can be found in Ref. \cite{ovseev0}.

This continuity implies, in particular, that the flow is defined uniquely everywhere, although the control $u(x)$ is defined uniquely only outside the hypersurface 
$\{\langle{B,{\partial\rho}/{\partial x}\rangle}{=}0\}$.
A similar phenomenon was discovered by I.A. Bogaevskii \cite{bogaev} for the gradient differential equations $\dot x=-{\partial f}/{\partial x}$, where $f$ is a nonsmooth convex function.

The singular part of the right-hand side of (\ref{sing_sys}) has the form $-\alpha(x)\frac{\partial f}{\partial x}$,
where $\alpha$ is a smooth nonnegative symmetric matrix, while  $f$  is a (nonsmooth) convex function.
The quadratic form
\begin{equation}
    \langle\alpha(x)\xi,\xi\rangle+\langle x,\xi\rangle^2
\end{equation}
is strictly positive, and the singular part of (\ref{sing_sys}) is invariant under scaling $x\mapsto\lambda x$ of the phase space.
These facts allows us  obtain differential inequalities for $\langle\alpha(x)\frac{\partial v}{\partial x},\frac{\partial v}{\partial x}\rangle$ and $\langle x,\frac{\partial v}{\partial x}\rangle^2$,
where $v$ is a solution of (\ref{transport}), and a priori bounds for the Lipschitz constant of $v$ in any domain of the form $\{|\phi_{t}(x)|\geq c\}$.
More precisely, we need differential inequalities for a solution of regularized equation (\ref{transport}), where $\sign$ is replaced by a smooth function approximating $\sign$ in $L_1$.
In this way, we obtain the Lipschitz constants which do not depend on the accuracy of the regularization.

\section*{Conclusion}

The generalized dry-friction control for a system of oscillators  is especially interesting because of its asymptotic optimality, previously studied in Refs. \cite{ovseev0,ovseev2,ovseev3}.
As it is usual in the optimal control theory, this control is not defined uniquely everywhere and is a discontinuous function of the phase state of the system.
Existence of the motion under control follows from the Filippov theory. However this theory does not answer  the question of uniqueness of trajectories.
It is shown that the existence and {\it uniqueness} issue of the motion under the generalized dry friction control can be resolved in the framework of DiPerna--Lions theory.
It is important and interesting to explore similar issues for the optimal control.

In real life, usage of a discontinuous control is accompanied with  difficulties of various kinds.
In the considered case, the discontinuous control in the form of a generalized dry friction arise naturally from the study of the asymptotic behavior of reachable sets of general linear systems.
As noted above, in the proof of the main Theorem \ref{motion2} we use an approximation of a singular law of motion by smoothing the function $x\mapsto\sign(x)$.
This quite common method, of course, can be useful in practical applications.

We emphasize that the control law in the generalized dry friction form does not completely solve the problem of damping the oscillations of a system of linear oscillators.
The method proposed in author's papers \cite{ovseev0,ovseev2,ovseev3} includes a combination of three control strategies:
for large, intermediate and small values of the energy of the system.
In the case of high energy, generalized dry friction is used. For intermediate values of the energy we use a similar control with a smaller amplitude.
By decreasing the amplitude of the control we are reducing the size of the standstill zones.
In a close vicinity of the equilibrium, this allows us  to switch to a completely different type of control based on common Lyapunov functions.

Thus, the result of the paper on the existence and uniqueness of the motion under the generalized dry friction is relevant at the first two stages of the motion to the terminal set,
while at the last stage the issue has a different nature because of the substantially different control law.
From the viewpoint of asymptotic optimality, the motion at the first stage plays a decisive role.

\section*{Acknowledgements}

This work is supported by the Russian Foundation for Basic Research (projects 14-08-00606 and 14-01-00476).

\end{document}